\documentclass[12pt]{article}

\usepackage{amssymb}
\usepackage{amsmath}

\def\shadowbox{\hbox{\rule[-0.0ex]{0.1ex}{1.2ex}%
\hspace{-0.1ex}\rule[-0.0ex]{1.2ex}{0.1ex}%
\hspace{0.0ex}\rule[-0.0ex]{0.1ex}{1.2ex}\hspace{-1.3ex}%
\rule[1.15ex]{1.25ex}{0.1ex}\hspace{-0.0ex}\rule[-0.25ex]{0.3ex}{1.1ex}%
\hspace{-1.2ex}\rule[-0.25ex]{1.1ex}{0.25ex}}}
\def\qed{\ifmmode \hbox{\hfill\shadowbox}
     \else \hphantom{x}\hfill\shadowbox \fi}

\newtheorem{theorem}{Theorem}[section]

\def\remark{{\noindent \bf Remark:\hspace{0.5em}}}

\def\Cst{{\mathbb C}}

\def\Rst{{\mathbb R}}
\def\Zst{{\mathbb Z}}

\def\theta{\vartheta}
\def\phi{\varphi}
\def\Cst{{\mathbb C}}

\def\Rst{{\mathbb R}}

\def\Zst{{\mathbb Z}}
\def\Lsp{{\boldsymbol L}}

\def\lsp{{\boldsymbol\ell}}
\def\Ltsp{{\Lsp^2}}
\def\Lpsp{{\Lsp^p}}
\def\Linsp{{\Lsp^{\infty}}}
\def\LtR{{\Lsp^2(\Rst)}}
\def\ltZ{{\lsp^2(\Zst)}}

\def\Res{\operatorname{Res}}
\def\sinc{\operatorname{sinc}}

\def\oh{{\frac{1}{2}}}

\newcommand{\gtight}{g^{t}}
\newcommand{\gxy}{g_{x,y}}
\newcommand{\gnm}{g_{na,mb}}

\newcommand{\Mg}{M_g}
\newcommand{\Mga}{M^{\ast}_g}
\newcommand{\tfnu}{\theta_4(\pi \nu;e^{-\pi \gamma})}
\newcommand{\tft}{\theta_4(\pi t;e^{-\pi/ \gamma})}
\renewcommand{\ggg}{g_{1,\gamma}}
\newcommand{\gtgg}{g^{t}_{1,\gamma}}
\newcommand{\gthg}{g^{t}_{2,\gamma}}

\newcommand{\gdhg}{g^{d}_{2,\gamma}}
\newcommand{\ghg}{g_{2,\gamma}}
\newcommand{\Zggg}{Z\ggg}
\newcommand{\Zghg}{Z\ghg}

\newcommand{\Zgdhg}{Z\gdhg}

\def\toinf{{\rightarrow \infty }}

\begin{document}

\title{Hyperbolic secants yield Gabor frames}

\author{A.J.E.M.~Janssen\thanks{Philips Research Laboratories WY-81, 5656 AA 
Eindhoven, The Netherlands; Email: a.j.e.m.janssen@philips.com} and Thomas 
Strohmer\thanks{Department of Mathematics, University of California, Davis, 
CA 95616-8633, USA; Email: strohmer@math.ucdavis.edu. T.S. acknowledges 
support from NSF grant 9973373.}}

\date{}
\maketitle



\begin{abstract}
We show that $(g_2,a,b)$ is a Gabor frame when $a>0, b>0, ab<1$ and
$g_2(t)=(\frac{1}{2}\pi \gamma)^{\frac{1}{2}} (\cosh \pi \gamma t)^{-1}$
is a hyperbolic secant with scaling parameter $\gamma >0$. This is
accomplished by expressing the Zak transform of $g_2$ in terms of
the Zak transform of the Gaussian
$g_1(t)=(2\gamma)^{\frac{1}{4}} \exp (-\pi \gamma t^2)$, together with
an appropriate use of the Ron-Shen criterion for being a Gabor frame. 
As a side result it follows that the windows, generating tight Gabor
frames, that are canonically associated to $g_2$ and $g_1$
are the same at critical density $a=b=1$. Also, we display
the ``singular'' dual function corresponding to the hyperbolic secant
at critical density.
\end{abstract}

\noindent
{\em AMS Subject Classification:} 42C15, 33D10, 94A12.

\noindent
{\em Key words:} Gabor frame, Zak transform, hyperbolic secant, theta
functions.

\section{Introduction and results}
\label{s:1}

Let $a>0, b>0$, and let $g \in \LtR$. A Gabor system $(g,a,b)$ consists
of all time- and- frequency shifted functions $\gnm$ with
integer $n, m$, where for $x,y \in \Rst$ we denote
\begin{equation}
\gxy(t) = e^{2\pi i yt} g(t-x),\,\,\, t\in \Rst.
\label{e1}
\end{equation}
We say that $(g,a,b)$ is a Gabor frame when there are $A>0, B< \infty$
(lower, upper frame bound, respectively) such that 
for all $f \in \LtR$ we have
\begin{equation}
A \|f\|^2 \le \sum_{n,m} |\langle f, \gnm \rangle|^2 \le B\|f\|^2.
\label{e2}
\end{equation}
Here $\| \cdot  \|$ and $\langle \cdot, \cdot\rangle$ denote the standard 
norm and inner product of $\LtR$. We refer to~\cite{Dau90,Dau92,FS98} for
generalities about frames for a Hilbert space and for both basic and
in-depth information about Gabor frames, frame operators, dual frames,
tight frames, etc. For a very recent and
comprehensive treatment of Gabor frames, we refer to~\cite{Gro01},
Chs.~5-9,~11-13; many of the more advanced results of modern Gabor theory
are covered there in a unified manner.

It is well known that a triple $(g,a,b)$ cannot be a Gabor frame when
$ab>1$, see for instance~\cite{Gro01}, Corollary~7.5.1. Also such a
system cannot be a Gabor frame when $ab=1$ and $g' \in \LtR, tg(t) \in
\LtR$ (extended Balian-Low theorem; see for instance~\cite{FS98}, Ch.~2).

In this paper we demonstrate that
$(g_2,a,b)$ is a Gabor frame when $ab<1$ and
\begin{equation}
g_{2,\gamma}(t)= \Big(\frac{\pi \gamma}{2}\Big)^{\frac{1}{2}} 
                 \frac{1}{\cosh \pi \gamma t}, \qquad t \in \Rst,
\label{e3}
\end{equation}
with $\gamma>0$ (normalization such that $\|g_{2,\gamma}\|=1$). The result
is obtained by (i) relating the Zak transform
\begin{equation}
(Zg)(t,\nu) = \sum_{l=-\infty}^{\infty} g(t-l) e^{2\pi i l \nu}, \qquad
t, \nu \in \Rst,
\label{e4}
\end{equation}
of $g=\ghg$ to the Zak transform of $g=\ggg$, where $\ggg$ is the
normalized Gaussian
\begin{equation}
\ggg (t)=(2\gamma)^{\frac{1}{4}} e^{-\pi \gamma t^2}, \qquad t \in \Rst,
\label{e5}
\end{equation}
(ii) using the observation that $(\ggg,a,b)$ is a Gabor frame and
(iii) an appropriate use of the Ron-Shen time domain 
criterion~\cite{RS97a} for a Gabor system $(g,a,b)$ to be a Gabor frame.
Explicitly, we show that there is a positive constant $E$ such that
\begin{equation}
(\Zghg)(t,\nu)=\frac{E (\Zggg)(t,\nu)}{\tfnu \tft} , \quad t,\nu \in \Rst.
\label{e6}
\end{equation}
Here $\theta_4(z;q)$ is the theta function, see~\cite{WW27}, Ch.~21,
\begin{equation}
\theta_4(z;q) = \sum_{n=-\infty}^{\infty} (-1)^n q^{n^2} e^{2 i nz},
\qquad z \in \Cst.
\label{e7}
\end{equation}
This $\theta_4$ is positive and bounded for real arguments $z$ and
parameter $q \in (0,1)$.

The rest of the paper is organized as follows.
In Section~\ref{s:2} we present the proof or our main result. In
Section~\ref{s:3} we prove formula~\eqref{e6} by using elementary 
properties of theta functions and basic complex
analysis. These same methods, combined with formula~\eqref{e6}, allow
us to compute explicitly the ``singular'' dual function $\gdhg$
(at critical density $a=b=1$) canonically associated to $\ghg$ according to
\begin{equation}
\gdhg(t)= \int \limits_{0}^{1} \frac{d\nu}{(\Zghg)^{\ast}(t,\nu)},
\qquad t \in \Rst. 
\label{e8}
\end{equation}
This is briefly presented in Section~\ref{s:4}. We also show in
Section~\ref{s:4} that the tight Gabor frame generating windows
$\gtgg, \gthg$, canonically associated to $\ggg, \ghg$ and given
at critical density $a=b=1$ by
\begin{equation}
\gtight(t)=\int \limits_{0}^{1} \frac{(Zg)(t,\nu)}{|(Zg(t,\nu)|} d\nu,
\qquad t \in \Rst,
\label{e9}
\end{equation}
for $g=\ggg, \ghg$, are the same. In~\cite{Jan01} a comprehensive
study is made of the process, embodied by formula~\eqref{e9}, of
passing from windows $g$ with few zeros in the Zak transform domain
to tight Gabor frame generating windows $\gtight$ at critical 
density. One of the observations in~\cite{Jan01} is that the
operation in~\eqref{e9} seems to diminish distances between positive, even,
unimodal windows enormously. The fact that $\gtgg=\gthg$ is an absurdly
accurate illustration of this phenomenon.

\section{Proof of the main result} \label{s:2}

In this section we present a proof for the result that $(\ghg,a,b)$ is
a Gabor frame when $a>0, b>0, ab<1$ and $\ghg$ is given by~\eqref{e3}.
According to the Ron-Shen criterion in the
time-domain~\cite{RS97a} we have that $(g,a,b)$ is a Gabor frame with
frame bounds $A>0, B<\infty$ if and only if
\begin{equation}
AI \le \frac{1}{b}\Mg(t) \Mga(t) \le B I, \qquad \text{a.e.~$t \in \Rst$}.
\label{e10}
\end{equation}
Here $I$ is the identity operator of $\ltZ$ and $\Mg(t)$ is the linear
operator of $\ltZ$ (in the notation of~\cite{FS98}, Subsec.~1.3.2),
whose matrix with respect to the standard basis of $\ltZ$ is given by
\begin{equation}
\Mg(t)=\big(g(t-na-l/b)\big)_{l \in \Zst, n \in \Zst}, \qquad
\text{a.e.~$t \in \Rst$}
\label{e11}
\end{equation}
(row index $l$, column index $n$). Since $\ghg$ is rapidly decaying, the
finite frame upper bound condition is easily seen to be satisfied, and we
therefore concentrate on the positive lower frame bound condition. We may
restrict here to the case where $a<1, b=1$, since $(g,a,b)$ is a Gabor
frame if and only if $(D_c g, a/c,bc)$ is a Gabor frame. Here $D_c$
is the dilation operator $(D_c f)(t)=c^{\frac{1}{2}} f(ct), t \in \Rst$,
defined for $f \in \LtR$ when $c>0$. Since $D_c \ghg = g_{2,\gamma c}$,
we only need to replace $\gamma>0$ by $\gamma b >0$ when $b \neq 1$. Hence
we shall show that there is an $A>0$ such that 
\begin{equation}
\sum_{n =-\infty}^{\infty}\Big|\sum_{l=-\infty}^{\infty}c_l g(t-na-l)\Big|^2 
\ge A \|\underline{c}\|^2, \quad \underline{c}=(c_l)_{l \in \Zst} \in \ltZ, 
t \in \Rst,
\label{e12}
\end{equation}
with $g=\ghg$.

Taking $\underline{c} \in \ltZ$ it follows from Parseval's theorem for
Fourier series and the definition of the Zak transform in~\eqref{e4} that
for any $n \in \Zst$
\begin{equation}
\sum_{l=-\infty}^{\infty} c_l g(t-na-l) =
\int \limits_{0}^{1}  (Zg)(t-na,\nu) C^{\ast}(\nu) d \nu, 
\label{e13}
\end{equation}
where $C(\nu)$ is defined by
\begin{equation}
C(\nu)= \sum_{l=-\infty}^{\infty} c^{\ast}_l e^{2\pi i l \nu},
\qquad \text{a.e.~$\nu \in \Rst$.}
\label{e14}
\end{equation}

Now assuming the result~\eqref{e6} (with $E>0$) we have that 
\begin{gather}
\sum_{l=-\infty}^{\infty} c_l \ghg(t-na-l) = 
\int \limits_{0}^{1} (\Zghg)(t-na,\nu) C^{\ast}(\nu) d\nu \notag \\
= \frac{E}{\theta_4(\pi(t-na);e^{-\pi/\gamma})} 
\int \limits_{0}^{1} (\Zggg)(t-na,\nu) \frac{C^{\ast}(\nu)}{\tfnu} d\nu.
\label{e15}
\end{gather}
Define the 1-periodic function $D(\nu)$ by
\begin{equation}
D(\nu)=\sum_{l=-\infty}^{\infty} d^{\ast}_l e^{2\pi i l \nu}:=
\frac{C(\nu)}{\tfnu}, \qquad \nu \in \Rst. 
\label{e16}
\end{equation}
It is well known that $(\ggg,a,1)$ is a Gabor frame, see for 
instance~\cite{Gro01}, Theorem~7.5.3. Accordingly, see~\eqref{e12}, there
is an $A_{1,\gamma}>0$ such that
\begin{equation}
\sum_{n=-\infty}^{\infty} \Big|\sum_{l=-\infty}^{\infty}
d_l \ggg(t-na-l)\Big|^2 \ge A_{1,\gamma} \| \underline{d} \|^2,
\quad \underline{d} \in \ltZ, t \in \Rst.
\label{e17}
\end{equation}
Letting
\begin{equation}
m_{\delta}:=\underset{z \in \Rst}{\min}\, 
\frac{1}{\theta_4^2(z;e^{-\pi \delta})} >0
\label{e18}
\end{equation}
for $\delta >0$, we then see from~\eqref{e13}--\eqref{e18} that
\begin{gather}
\sum_{n=-\infty}^{\infty} 
\Big|\sum_{l=-\infty}^{\infty} c_l \ghg(t-na-l)\Big|^2 \ge \notag \\
m_{1/\gamma} E^2 \sum_{n=-\infty}^{\infty} 
\Big|\int \limits_{0}^{1} (\Zggg)(t-na,\nu) D^{\ast}(\nu) d\nu \Big|^2 =
\notag \\
m_{1/\gamma} E^2 \sum_{n=-\infty}^{\infty} \Big| \sum_{l=-\infty}^{\infty}
d_l \ggg(t-na-l)\Big|^2 \ge 
m_{1/\gamma} E^2 A_{1,\gamma} \| \underline{d} \|^2.
\label{e19}
\end{gather}
Finally, with $D$ and $C$ related as in~\eqref{e16} we have from
Parseval's theorem for Fourier series that
\begin{gather}
\| \underline{d} \|^2 = \int \limits_{0}^{1} |D(\nu)|^2 d\nu =
\int \limits_{0}^{1} \frac{1}{|\tfnu|^2} |C(\nu)|^2 d\nu \notag \\
\ge m_{\gamma}\int \limits_{0}^{1} |C(\nu)|^2 d\nu = m_{\gamma} 
\|\underline{c}\|^2.  
\label{20}
\end{gather}
Hence
\begin{equation}
\sum_{n=-\infty}^{\infty} 
|\sum_{l=-\infty}^{\infty} c_l \ghg(t-na-l)|^2 \ge 
m_{\gamma} m_{1/\gamma} E^2 A_{1,\gamma} \|\underline{c}\|^2,
\label{e21}
\end{equation}
as required.

\section{Expressing $\Zghg$ in terms of $\Zggg$} \label{s:3}

In this section we express $Z \ghg$, with $\ghg$ given in~\eqref{e3}, 
in terms of $Z \ggg$, with $\ggg$ given in~\eqref{e5}. For this we need
some basic facts of the theory of theta functions as they can be found
in~\cite{WW27}, Ch.~21. The four theta functions are for $q \in (0,1)$
given by 
\begin{align}
\theta_1(z;q)=&\frac{1}{i} \sum_{n=-\infty}^{\infty} (-1)^n q^{(n+\oh)^2}
e^{(2n+1)iz}, \label{e22} \\
\theta_2(z;q)=&\sum_{n=-\infty}^{\infty} q^{(n+\oh)^2}
e^{(2n+1)iz}, \label{e23} \\
\theta_3(z;q)=&\sum_{n=-\infty}^{\infty} q^{n^2} e^{2niz}, \label{e24} \\
\theta_4(z;q)=&\sum_{n=-\infty}^{\infty} (-1)^n q^{n^2} e^{2niz}, \label{e25} 
\end{align}
where $z \in \Cst$. One thus easily gets from~\eqref{e4} and~\eqref{e5}
that
\begin{equation}
(\Zggg)(t,\nu)=(2\gamma)^{\frac{1}{4}} e^{-\pi \gamma t^2}
\theta_3(\pi(\nu-i\gamma t);e^{-\pi \gamma}), \quad t,\nu \in \Rst.
\label{e26}
\end{equation}

\begin{theorem}
\label{th1}
We have for $t,\nu \in \Rst$
\begin{align}
(\Zghg)(t,\nu) &= 2^{-\oh} \pi^{\oh} \gamma \theta'_1(0;e^{-\pi \gamma})
e^{-\pi \gamma t^2} 
\frac{\theta_3 (\pi (\nu-i\gamma t);e^{-\pi \gamma})}{\tfnu \tft} \notag \\
&= \pi^{\oh} \Big(\frac{\gamma}{2}\Big)^{\frac{3}{4}} 
\theta'_1(0;e^{-\pi \gamma}) \frac{(\Zggg)(t,\nu)}{\tfnu \tft}. 
\label{e27}
\end{align}
\end{theorem}
\begin{proof}
For fixed $t \in \Rst$ we compute the Fourier coefficients $b_n(t)$
of the 1-periodic function
\begin{equation}
\frac{\theta_3 (\pi (\nu-i\gamma t);e^{-\pi \gamma})}{\tfnu} 
=\sum_{n=-\infty}^{\infty} b_n(t) e^{2\pi i n \nu}, \quad \nu \in \Rst,
\label{e28}
\end{equation}
by the method that can be found in~\cite{WW27}, Sec.~22.6. For brevity
we shall suppress the expression ``$;e^{-\pi \gamma}$'' in 
$\theta_3(\pi(\nu-i\gamma t);e^{-\pi \gamma})$, etc.\ in the
remainder of this proof.

We thus have
\begin{equation}
b_n(t)=\int \limits_{-\oh}^{\oh} 
\frac{\theta_3(\pi(\nu-i\gamma t))}{\theta_4(\pi \nu)} e^{-2\pi i n \nu}
d\nu, \quad n \in \Zst,  
\label{e29}
\end{equation}
and we consider for $n \in \Zst$
\begin{equation}
c_n(t)=\int \limits_{C} 
\frac{\theta_3(\pi(\nu-i\gamma t))}{\theta_4(\pi \nu)} e^{-2\pi i n \nu}
d\nu,
\label{e30}
\end{equation}
where $C$ is the edge of the rectangle with corner points $-\oh,\oh,
\oh+i\gamma, -\oh+i\gamma$, taken with positive orientation. It follows
from~\cite{WW27}, Sec.~21.12 that the function
$\theta_4(\pi \nu)$ has a first order zero at $\nu=\oh i \gamma$ and no zeros
elsewhere on or within $C$. Therefore, by Cauchy's theorem,
\begin{equation}
c_n(t) = 2\pi i \underset{\nu = \oh i \gamma}{\Res} \left[
\frac{\theta_3(\pi(\nu-i\gamma t))}{\theta_4(\pi \nu)}e^{-2\pi i n \nu}
\right] =
\frac{2i \theta_3(\pi i(\oh - t)\gamma)e^{\pi \gamma n}}
     {\theta_4'(\oh \pi i \gamma)}.
\label{e31}
\end{equation}
By~\cite{WW27}, Ex.~2, first identity, on p.~464 we have 
($\tau=i\gamma,q=\exp(-\pi \gamma),z=0)$
\begin{equation}
\theta'_4\Big(\oh \pi i \gamma\Big)=i e^{\frac{1}{4}\pi \gamma} \theta'_1(0).
\label{e32}
\end{equation}
Next, by the first formula on p.~475 of~\cite{WW27}
($z=\pi i (\oh-t)\gamma, \tau'=-1/\tau=i/\gamma, q'=\exp(\pi i \tau')=
\exp(-\pi/\gamma)$)
\begin{equation}
\theta_3\big(\pi i(\oh -t)\gamma\big)=\gamma^{-\oh} 
\exp\big(\pi \gamma (\oh-t)^2\big) 
\theta_3\big(\pi (\oh-t);e^{-\pi/\gamma}\big).
\label{e33}
\end{equation}
and by~\cite{WW27}, Ex.~2, fourth identity, on p.~464
($z=-\pi t$, $q'$ instead of $q$)
\begin{equation}
\theta_3\big(\pi (\oh -t);e^{-\pi/\gamma}\big)=
\theta_4(-\pi t;e^{-\pi/\gamma})= \theta_4(\pi t;e^{-\pi/\gamma}),
\label{e34}
\end{equation}
where we also have used that $\theta_4$ is an even function. It thus
follows from~\eqref{e31}--\eqref{e34} that
\begin{equation}
c_n(t)=\frac{2}{\theta'_1(0)} \gamma^{-\oh} e^{-\frac{1}{4}\pi \gamma}
e^{\pi \gamma (\oh-t)^2} \tft e^{\pi \gamma n}.
\label{e35}
\end{equation}

On the other hand, we have by 1-periodicity of the integrand in~\eqref{e30}
(in $\nu$) that the two integrals along the vertical edges of $C$ cancel
one another. Hence
\begin{gather}
c_n(t)=\int \limits_{-\oh}^{\oh} 
\frac{\theta_3(\pi(\nu-i\gamma t))}{\theta_4(\pi \nu)} 
e^{-2\pi i n \nu}d\nu - \int \limits_{-\oh+i\gamma}^{\oh+i\gamma} 
\frac{\theta_3(\pi(\nu-i\gamma t))}{\theta_4(\pi \nu)} 
e^{-2\pi i n \nu}d\nu  \notag \\
= b_n(t) - e^{2\pi \gamma n} \int \limits_{-\oh}^{\oh}
\frac{\theta_3(\pi(\nu-i\gamma t)+\pi i\gamma)}{\theta_4(\pi\nu+\pi i\gamma)} 
e^{-2\pi i n \nu}d\nu.
\label{e36}
\end{gather}
Furthermore, by the table in~\cite{WW27}, Ex.~3 on p.~465 we have 
\begin{equation}
\theta_3\big(\pi (\nu-i\gamma t)+\pi i \gamma\big)=
e^{\pi \gamma} e^{-2i \pi (\nu-i\gamma t)} \theta_3 \big(\pi (\nu-i \gamma t)\big),
\label{e37}
\end{equation}
and
\begin{equation}
\theta_4(\pi \nu +\pi i \gamma) = - e^{\pi \gamma} e ^{-2 i \pi \nu}
\theta_4(\pi \nu).
\label{e38}
\end{equation}
It thus follows that
\begin{gather}
c_n(t)= b_n(t) + e^{-2\pi \gamma (t-n)} \int \limits_{-\oh}^{\oh}
\frac{\theta_3(\pi(\nu-i\gamma t))}{\theta_4(\pi\nu)} 
e^{-2\pi i n \nu}d\nu \notag \\
= b_n(t) (1+e^{-2\pi \gamma (t-n)}).
\label{e39}
\end{gather}
By~\eqref{e35} we then find that
\begin{align}
b_n(t)= &\frac{2\gamma^{-\oh}}{\theta'_1(0)} e^{-\frac{1}{4}\pi \gamma}
e^{\pi \gamma(\oh-t)^2} \tft \frac{e^{\pi \gamma n}}{1+e^{-2\pi \gamma(t-n)}}
\notag \\
=&\frac{\gamma^{-\oh}}{\theta'_1(0)} e^{-\frac{1}{4}\pi \gamma}
e^{\pi \gamma (\oh-t)^2+\pi \gamma t} \frac{\tft}{\cosh \pi \gamma (t-n)}
\notag \\
=&\frac{\gamma^{-\oh}}{\theta'_1(0)} e^{\pi \gamma t^2}
\frac{\tft}{\cosh \pi \gamma (t-n)}.
\label{e40}
\end{align}
We thus conclude, see~\eqref{e28}, that
\begin{equation}
\frac{\theta_3(\pi (\nu-i\gamma t))}{\theta_4(\pi \nu)}
=\frac{\gamma^{-\oh}}{\theta'_1(0)} e^{\pi \gamma t^2} \tft 
\sum_{n=-\infty}^{\infty}\frac{e^{2\pi i n \nu}}{\cosh \pi \gamma (t-n)}.
\label{e41}
\end{equation}
That is,
\begin{align}
(\Zghg)(t,\nu) &=
\Big(\frac{\pi \gamma}{2}\Big)^{\oh}
\sum_{n=-\infty}^{\infty}\frac{e^{2\pi i n \nu}}{\cosh \pi \gamma (t-n)}
\notag \\
&= 2^{-\oh} \pi^{\oh} \gamma \theta'_1(0)
\frac{e^{-\pi \gamma t^2} \theta_3(\pi(\nu-i\gamma t))}
{\theta_4(\pi \nu) \tft},
\label{e42}
\end{align}
and this is the first line identity in~\eqref{e27}. The second line
identity in~\eqref{e27} follows from~\eqref{e26}. This completes the proof.
\end{proof}

\underline{\bf Note:} We have that
\begin{equation}
\theta'_1(0;e^{-\pi \gamma}), \tfnu, \tft >0
\label{e43}
\end{equation}
for $\nu \in \Rst, t \in \Rst$. Indeed, we have by~\cite{WW27}, Sec.21.41
that
\begin{equation}
\theta'_1(0;e^{-\pi \gamma})= \theta_2(0;e^{-\pi \gamma})
\theta_3(0;e^{-\pi \gamma}) \theta_4(0;e^{-\pi \gamma}).
\label{e44} 
\end{equation}
One sees directly from~\eqref{e23} that $\theta_2(0;q)>0$. Also,
from the formula on p.~476 of~\cite{WW27} just before Ex.~1, one sees that
\begin{equation}
\theta_3(z;q)=\theta_4\Big(z+\oh \pi;q\Big)>0,\qquad z \in \Rst.
\label{e45}
\end{equation}

\remark Theorem~\ref{th1} implies that
$(Z g_2)(t,\nu)/(Z g_1)(t,\nu)$ can be factored into a function of
$t$ and a function of $\nu$. This factorization is crucial in the proof
of the main result in Section~\ref{s:2}. It seems possible to extend the
approach in Section~\ref{s:2} to other pairs of windows $g_1, g_2$
where $(g_1,a,b)$ is a frame and $Z g_2/Z g_1$ (nearly) factorizes.
We do not pursue this extension in this paper.

\section{Canonical dual window and tight window at critical density} 
\label{s:4}

Theorem~\ref{th1} allows us to calculate the Zak transform of the
canonical dual $\gdhg$ of $\ghg$ for rational values of $ab<1$
using Zibulski-Zeevi matrices representing the
frame operator in the Zak transform domain, see for instance~\cite{Jan98},
Sec.~1.5, \cite{Gro01}, Sec.~8.3, and, of course, \cite{ZZ93}. The dual
window can then be obtained as
\begin{equation}
\gdhg(t)= \int \limits_{0}^{1}  (\Zgdhg)(t,\nu) d\nu, \qquad t \in \Rst.
\label{e46}
\end{equation}
For $a=b=1$ the Zibulski-Zeevi matrices are just scalars, viz.\
$|(\Zghg)(t,\nu)|^2$, and we have formally
\begin{equation}
(\Zgdhg)(t,\nu)= \frac{1}{(\Zghg)^{\ast}(t,\nu)}, \qquad t,\nu \in \Rst,
\label{e47}
\end{equation}
and
\begin{equation}
\gdhg(t)= \int \limits_{0}^{1}  \frac{d\nu}{(\Zghg)^{\ast}(t,\nu)}, 
\qquad t \in \Rst.
\label{e48}
\end{equation}
The identities here hold only formally since $1/\Zghg \neq
\Ltsp_{\text{loc}}(\Rst^2)$, and, as said, $(\ghg,1,1)$ is not a frame.
However, nothing prevents us from computing the right-hand side
of~\eqref{e48} for $t \in \Rst$ not of the form $n+\oh, n \in \Zst$. 
For these $t$ we have that $(\Zghg)(t,\nu)$ is a well-behaved zero-free
function of $\nu \in \Rst$.

When we do this calculation for the case that $\gamma=1$, we find,
using the same methods as in the proof of Theorem~\ref{th1}, that for
$t \in (-\frac{1}{2},\frac{1}{2}) , n \in \Zst$
\begin{equation}
g^{d}_{2,1}(t+n)=
\frac{2^{\oh} (-1)^n \theta_4^2(\pi t) e^{2\pi t^2 +2\pi n t}}
{\pi^{\oh}(\theta'_1(0))^2\cosh \pi (t+n)}\, .
\label{e49}
\end{equation}
Here all theta functions are with parameter $q=e^{-\pi}$. The 
same procedure for the Gaussian $g_{1,1}$ yields for 
$t \in \Rst$ not of the form $n+1/2, n \in \Zst$
\begin{equation}
g^{d}_{1,1}(t)=\frac{2^{-\frac{1}{4}}}{\theta'_1(0)} e^{\pi t^2}
\sum_{n-\oh \ge |t|} (-1)^n e^{-\pi (n-\oh)^2},
\label{e50}
\end{equation}
see for instance~\cite{Bas80a}, \cite{Jan81a}, Subsec.~2.14, 
and~\cite{Jan82}, Subsec.~4.4. Although the two functions 
in~\eqref{e49}--\eqref{e50} seem quite different, one can show that both
functions are bounded but not in $\Lpsp(\Rst)$ when $1\le p < \infty$.
Moreover, there is for both functions exponential decay as $|t|\toinf$ 
away from the set of half-integers.

In a similar fashion one can consider the canonical tight windows 
$\gtight$ associated with $g=\ggg, \ghg$ at critical density $a=b=1$. We
have that these windows are given in the Zak transform domain through the
formula $Z \gtight = Zg/|Zg|$, whence
\begin{equation}
\gtight(t)=\int \limits_{0}^{1} \frac{(Zg)(t,\nu)}{|(Zg)(t,\nu)|} d\nu,
\qquad t \in \Rst. 
\label{e51}
\end{equation}
In this case the definitions are more than just formal since
$Zg/|Zg| \in \Linsp(\Rst^2)$, while one can show, see~\cite{Jan01}
for details, that in $\Ltsp$-sense
\begin{equation}
\gtight = \underset{a \uparrow 1}{\lim}\, S^{-\oh}_a g,
\label{e52}
\end{equation}
where for $g = \ggg, \ghg$ we have denoted the frame operator corresponding
to the Gabor frame $(g,a,a)$ by $S_a$. Now as a surprising consequence
of Theorem~\ref{th1} we have that $\gtgg = \gthg$ since 
$\Zggg/\Zghg$ is positive everywhere. In~\cite{Jan01} it has
been observed that the mapping $g \rightarrow \gtight$ in~\eqref{e51}
diminishes distances between even, positive, unimodal windows $g$
enormously when scaling parameters are set correctly. The fact that
$\gtgg=\gthg$ is an extremely accurate demonstration of this phenomenon.

The systems $(\gtight,1,1)$ form orthonormal bases for $\LtR$, see,
for instance, \cite{Gro01}, Corollary~7.5.2. Interestingly, in the
case that $g_{\gamma} = \ggg$ or $\ghg$, we have that, in $\Ltsp$-sense,
\begin{equation}
\underset{\gamma \downarrow 0}{\lim}\, \gtight_{\gamma} = \sinc \pi \cdot ,
\qquad 
\underset{\gamma \rightarrow \infty}{\lim} 
\,\gtight_{\gamma} = \chi_{(-\oh,\oh)}.
\label{e53}
\end{equation}
Hence we have a family $\gtight_{\gamma}, \gamma>0$ of reasonably behaved
tight frame generating windows at critical density $a=b=1$ that
interpolates between the Haar window $\chi_{(-\oh,\oh)}$ and its
Fourier transform $\sinc \pi \cdot$ as $\gamma$ varies between
$\infty$ and $0$. See~\cite{Jan01} for details.


\end{document}